\documentclass[12pt]{article}
\usepackage{mathrsfs}
\usepackage{amsfonts}
\usepackage{amssymb}
\usepackage{amsthm}
\usepackage{enumerate}
\usepackage{latexsym,amsmath}
\usepackage{graphicx}
\usepackage{abstract}
\usepackage{amsmath}
\usepackage{caption}
\title{On the generalized distance spectral radius of graphs}

\date{}
\begin{document}
\author{Shu-Yu Cui$^{a,b}$, Gui-Xian Tian$^b$\footnote{Corresponding author. E-mail: gxtian@zjnu.cn or guixiantian@163.com(G.-X. Tian)}, Lu Zheng$^b$\\
       {\small{\it $^a$Xingzhi College, Zhejiang Normal University, Jinhua, 321004, P.R. China}}\\
       {\small{\it $^b$ Department of Mathematics, Zhejiang Normal University, Jinhua, 321004,P.R. China}}
}\maketitle
\begin{abstract}

The generalized distance spectral radius of a connected graph $G$ is
the spectral radius of the generalized distance matrix of $G$,
defined by
$$D_\alpha(G)=\alpha Tr(G)+(1-\alpha)D(G), \;\;0 \le
\alpha \le 1,$$ where $D(G)$ and $Tr(G)$ denote the distance matrix
and diagonal matrix of the vertex transmissions of $G$,
respectively. This paper characterizes the unique graph with minimum
generalized distance spectral radius among the connected graphs with
fixed chromatic number, which answers a question about the
generalized distance spectral radius in spectral extremal theories.
In addition, we also determine graphs with minimum generalized
distance spectral radius among the $n$-vertex trees and unicyclic
graphs, respectively. These results generalize some known results
about distance spectral radius and distance signless Laplacian
spectral radius of graphs.


\end{abstract}

\section{Introduction}

All graphs considered are finite, simple and connected in this
paper. Let $G$ be a connected graph with vertex set $V(G)$ and edge
set $E(G)$. The distance between two vertices $u$ and $v$ in $G$,
denoted by $d_{uv}$, is defined to the length of the shortest path
between $u$ and $v$. The transmission $Tr(u)$ of a vertex $u$ in $G$
is defined to be the sum of the distances from $u$ to all other
vertices in $G$, i.e, $Tr(u) = \sum\nolimits_{v \in V(G)} {d_{uv}}
$. We say that $G$ is $k$-transmission regular graph if $Tr(u)=k$
for any $u\in V(G)$. The distance matrix $D(G)=(d_{uv})$ of $G$ is
the matrix indexed by vertices of $G$ with its ($u,v$)-entry equal
to $d_{uv}$. Let $Tr(G)$ be the diagonal matrix of the vertex
transmissions in $G$. Then the distance signless Laplacian matrix
and distance Laplacian matrix are defined by Aouchiche and Hansen in
\cite{Aouchiche4} as $D^Q(G)=Tr(G)+D(G)$ and $D^L(G)=Tr(G)-D(G)$,
respectively. Recently, Cui et al.\cite{Cui2019} proposed to study
the convex linear combinations $D_\alpha(G)$ of $Tr(G)$ and $D(G)$,
defined by
\begin{equation}\label{1}
D_\alpha(G)=\alpha Tr(G)+(1-\alpha)D(G), \;\;0 \le \alpha  \le 1.
\end{equation}
Here $D_\alpha(G)$ is called the generalized distance matrix of $G$.
Observe that $ D_0(G) = D(G),\;\;D_\frac{{{1}}}{{{2}}}(G) =
\frac{1}{2}D^Q(G),\;\;D_1(G) =Tr(G),$ and $D_\alpha (G) - D_\beta
(G) = (\alpha  - \beta )D^L(G),$ which validates the rationality of
the definition (\ref{1}). In this setup, we can see the matrices
$D(G)$, $D^Q(G)$, $D^L(G)$ and $Tr(G)$ from a fresh perspective.

For an $n\times n$ matrix $M$ associated to a graph $G$, the
eigenvalues of $M$ are the zeros of the characteristic polynomial
$\det(xI_n-M)$ of $M$. The multiset of all eigenvalues of $M$ is
called the spectrum of $M$. In particular, the eigenvalues of $D(G)$
($D^Q(G)$ and $D_\alpha(G)$, respectively) are called the distance
(distance signless Laplacian and generalized distance, respectively)
eigenvalues of $G$, which are all real. The largest distance
(distance signless Laplacian and generalized distance, respectively)
eigenvalue of $G$ is called distance (distance signless Laplacian
and generalized distance, respectively) spectral radius of $G$,
denoted by $\rho_D(G)$ ($\rho_{D^Q}(G)$ and $\rho_{D_\alpha}(G)$,
respectively).

The distance spectrum of a connected graph has been investigated
extensively. For example, Stevanovi\'{c} and
Ili\'{c}\cite{Stevanovic2010} proved that the star $S_n$ is the
unique graph with minimum distance spectral radius among $n$-vertex
trees. Yu et al.\cite{Yu2011} determined the unique graphs with
maximum and minimum distance spectral radii among unicyclic graphs
with fixed number of vertices. Liu\cite{Liu2010} determined the
unique graphs with minimum distance spectral radius in three classes
of simple connected graphs with fixed vertex connectivity, matching
number and chromatic number, respectively. For more results on the
distance spectral radius, see the recent survey \cite{Aouchiche2014}
as well as the references therein. Recently, the distance Laplacian
spectrum and distance signless Laplacian spectrum of graphs have
also been studied in many papers. For example, Aouchiche and Hansen
\cite{Aouchiche5} showed that the star $S_n$ is the unique tree with
minimum distance Laplacian spectral radius. Xing et
al.\cite{Xing2014} determined the unique graphs with minimum
distance signless Laplacian spectral radius among the $n$-vertex
trees, unicyclic graphs and bipartite graphs, respectively. Li et
al.\cite{Li2014} gave a lower bound on the distance signless
Laplacian spectral radius in terms of chromatic number and
determined the unique graph with minimum distance signless Laplacian
spectral radius among connected graphs with fixed number of vertices
and chromatic number. For more review about distance (signless)
Laplacian spectral radius of graphs, readers may refer to
\cite{Aouchiche4,Aouchiche5,Das2011,Lin2016,Tian2015,Xing2014} and
the references therein.

Compared with the much studied distance (distance Laplacian and
distance signless Laplacian) spectra, the study of generalized distance spectrum has just
been proposed by Cui et al. in \cite{Cui2019}. Some basic spectral properties of generalized distance
matrix of graphs are established and bounds for the generalized distance spectral radius were obtained.
In \cite{Cui2019}, Cui et al. also determined the unique graph with minimum generalized distance spectral radius among
all connected bipartite graphs with fixed number of vertices. The following problem was proposed:
\emph{Characterize graphs with minimum generalized distance spectral radius among all connected $r$-partite graphs of order $n$
with $r\geq3$.} Recall that the chromatic number of a connected graph $G$ is the smallest number of colors
needed to color the vertices of $G$ such that two adjacent vertices have different colors. The vertex set in the same color class forms an
independent set of $G$. Thus, above problem is equivalent to the following problem:
\emph{Characterize graphs with minimum generalized distance spectral radius among all
connected graphs of order $n$ with fixed chromatic number $r\geq3$.}

In this paper, we give a lower bound on generalized distance spectral radius of connected graphs in terms of chromatic number
and determine the unique graph with minimum generalized distance spectral radius
among all connected graphs of order $n$ with fixed chromatic number. The result gives a positive answer on the proposed problem.
We also determine the unique graphs with minimum generalized distance spectral radius
among the $n$-vertex trees and unicyclic graphs, respectively. These results  generalize some known results
about distance spectral radius and distance signless Laplacian spectral radius of graphs.

\section{Minimum generalized distance spectral radius of graphs with fixed chromatic number}

Let $G$ be a connected graph with $V(G)=\{v_1,v_2,\ldots,v_n\}$ and
the transmission sequence $(Tr(v_1),Tr(v_2),\dots,Tr(v_n))$. For a
real vector $\textbf{x}=(x_1,x_2,\dots,x_n)^T$, the quadratic form
$\textbf{x}^TD_\alpha(G)\textbf{x}$ can be written as
\begin{equation}\label{2}
\textbf{x}^TD_\alpha(G)\textbf{x}=\alpha \sum\limits_{i=1}^n
{Tr(v_i)x^2_i + 2(1 - \alpha )} \sum\limits_{1 \le i < j \le n}
{d_{v_iv_j}x_ix_j},
\end{equation}
In addition, if $G$ is connected, then $D_\alpha(G)$ is a
nonnegative and irreducible matrix for any $0\leq \alpha<1$. By the
Perron Frobenius Theorem, $\rho_{D_\alpha}(G)$ has multiplicity one
and there exists a unique positive unit eigenvector
$\textbf{x}=(x_1,x_2,\dots,x_n)^T$ corresponding to it. We refer to
$\textbf{x}$ as to the Perron vector of $G$. Then the system of
eigenequations $D_\alpha(G)\textbf{x}=\rho_{D_\alpha}(G)\textbf{x}$
can be written as, for each vertex $v_i\in V(G)$,
    \begin{equation}\label{3}
    \rho_{D_\alpha}(G) x_i = \alpha Tr(v_i)x_i + (1 - \alpha )\sum\limits_{j=1}^n
    {d_{v_iv_j}x_j}.
    \end{equation}

For two nonadjacent vertices $u,v\in V(G)$, let $G+uv$ be the graph obtained from $G$ by adding the edge $uv$.
By the Perron Frobenius Theorem, we have the following lemma immediately.

\paragraph{Lemma 2.1.}(\cite{Cui2019,Minc1988}) {\it Let $G$ be a connected graph with two nonadjacent vertices $u,v\in V(G)$.
Then $\rho_{D_\alpha}(G+uv)<\rho_{D_\alpha}(G)$ for $0\leq \alpha<1$.}\\

Recall that two vertices $u$ and $v$ are equivalent in $G$ if
there exists an automorphism $\sigma:G \to G$ such that
$\sigma(u)=v$.

\paragraph{Lemma 2.2.}(\cite{Cui2019}) {\it Let $G$ be a connected graph of order
$n$ and $\textbf{x}=(x_1,x_2,\dots,x_n)^T$ be the Perron vector
corresponding to $\rho_{D_\alpha}(G)$. If two vertices $u$ and
$v$ are equivalent in $G$, then $x_u=x_v$.}

\paragraph{Lemma 2.3.} (\cite{Minc1988}) {\it Let $M$ be a nonnegative irreducible
 matrix of order $n$ with the spectral radius $\rho_M$ and row sum
$r_1(M),r_2(M),\ldots,r_n(M)$. Then
$${{\min\limits_{1\leq i\leq n}\{ r_i(M)\} }} \le \rho_M\le \max \limits_{1\leq i\leq n}\{ r_i(M)\}.$$
Moveover, either of the equalities holds if and only if the row sums of $M$ are all equal.}\\

Let $K_{n_1,n_2,\ldots,n_r}$ denote the complete $r$-partite graph
whose parts have size $n_1,n_2,\ldots,n_r$, respectively. Recall
that the Tur\'{a}n graph $T_{n,r}$ is a complete $r$-partite graph
on $n$ vertices for which the numbers of vertices of every parts are
as equal as possible.

Next we give a lower bound on generalized distance spectral radius
of connected graphs in terms of chromatic number and determine the
unique graph with minimum generalized distance spectral radius among
all connected graphs of order $n$ with fixed chromatic number. The
result answers the proposed problem in Introduction (also see
\cite{Cui2019}).

\paragraph{Theorem 2.4.} {\it Let $G$ be a connected graph of order $n$
with given chromatic number $r$, where $3 \leq r \leq n-1$. If
$0\leq\alpha\leq 1-\frac{1}{r}$, then
\begin{equation}\label{4}\small
\rho _{D_\alpha  } (G) \ge \rho _{D_\alpha  } (T_{n,r})=\frac{{n +
2d + \alpha n - 3 + \sqrt {(n(1 - \alpha ) - 1)^2  + 4s(1 - \alpha
)(d + 1)} }}{2},
\end{equation}
where $n=rd+s$, $0\leq s<r$. Furthermore, the equality in (\ref{4})
holds if and only if $G\cong T_{n,r}$.}
\begin{proof}
We first prove that $\rho _{D_\alpha  } (G)\geq\rho _{D_\alpha  }
(T_{n,r})$ with equality if and only if $G$ is the Tur\'{a}n graph
$T_{n,r}$. Let $G$ be the graph with minimum generalized distance
spectral radius among all connected graphs of order $n$ with
chromatic number $r$. Lemma 2.1 implies that $G$ is the complete
$r$-partite graph $K_{n_1,n_2,\ldots,n_r}$, where
$n=n_1+n_2+\cdots+n_r$. Without loss of generality, we assume that
$n_1\geq n_2\geq\cdots\geq n_r$.

Suppose towards contradiction that $G$ is not the Tur\'{a}n graph.
Then $n_1-n_r>1$. Consider the graph
$G'=K_{n_1-1,n_2,\ldots,n_{r-1},n_r+1}$ where $V_i'\;\;
(i=1,2,\ldots,r)$ is the vertex partition of the vertex set of $G'$
with $|V_1'|=n_1-1,\;\;|V_r'|=n_r+1,\;\;|V_i'|=n_i$ for
$i=2,\ldots,r-1$. Let $\textbf{x}$ be the unit Perron vector of $G'$
corresponding to $\rho_{D_\alpha}(G')$. In the light of Lemma 2.2,
the equivalent vertices have same entries in $\textbf{x}$. Thus
$\textbf{x}$ can be shaped as
$$\textbf{x}=(\underbrace {x_1,\ldots,x_1}_{n_1-1},\underbrace {x_2,\ldots
,x_2}_{n_2},\ldots,\underbrace {x_{r-1},\ldots
,x_{r-1}}_{n_{r-1}},\underbrace {x_r,\ldots,x_r}_{n_r+1}).$$

\paragraph{Claim 1.} {\it $x_1\geq x_r$.}\\
\\
\emph{Proof of Claim 1.} First, the identity (\ref{3}) implies that,
for any $u\in V_1'$ and $v\in V_r'$,
\[
\rho _{D_\alpha  } (G')x_1  = \alpha Tr(u)x_1  + (1 - \alpha )(2(n_1
- 2)x_1  + \sum\limits_{i = 2}^{r - 1} {n_i x_i }  + (n_r  + 1)x_r
),
\]
\[
\rho _{D_\alpha  } (G')x_r  = \alpha Tr(v)x_r  + (1 - \alpha )((n_1
- 1)x_1  + \sum\limits_{i = 2}^{r - 1} {n_i x_i }  + 2n_r x_r ).
\]
Combining with above two equation, one has
\[
(\rho _{D_\alpha  } (G') - \alpha Tr(u) - (1 - \alpha )(n_1  -
3))x_1  = (\rho _{D_\alpha  } (G') - \alpha Tr(v) - (1 - \alpha
)(n_r  - 1))x_r.
\]
Observe that $Tr(u)=n+n_1-3$ and $Tr(v)=n+n_r-1$. Then
\begin{equation}\label{5}
(\rho _{D_\alpha  } (G') - \alpha n - (n_1  - 3))x_1  = (\rho
_{D_\alpha  } (G') - \alpha n - (n_r  - 1))x_r.
\end{equation}
From Lemma 2.3, we gets
\[
\rho _{D_\alpha  } (G') \ge \min \{ n + n_1  - 3,n + n_{r - 1}  -
2,n + n_r  - 1\}  \ge n - 1.
\]
Thus
\begin{equation}\label{6}
\rho _{D_\alpha  } (G') - \alpha n - (n_r  - 1)\geq
(1-\alpha)n-n_r>0,
\end{equation}
where the last inequality holds as
$$\alpha\leq 1-\frac{1}{r}<1-\frac{n_r}{n_1+n_2+\cdots+n_r}=1-\frac{n_r}{n}.$$
It is easy to see that
\begin{equation}\label{6i}
\rho _{D_\alpha  } (G') - \alpha n - (n_1 - 3) \leq\rho _{D_\alpha
} (G') - \alpha n - (n_r  - 1).
\end{equation}
Therefore, from (\ref{5}), (\ref{6}) and (\ref{6i}), implying that
$x_1\geq x_r$. $\Box$

Since the respective transmission degree sequences of $G'$ and $G$
are
\[
\left( {(n + n_1  - 3)^{(n_1  - 1)} ,(n + n_2  - 2)^{(n_2 )} ,
\ldots ,(n + n_{r - 1}  - 2)^{(n_{r - 1} )} ,(n + n_r  - 1)^{(n_r  +
1)} } \right)
\]
and
\[
\left( {(n + n_1  - 2)^{(n_1  - 1)} ,(n + n_2  - 2)^{(n_2 )} ,
\ldots ,(n + n_r  - 2)^{(n_r )} ,(n + n_1  - 2)} \right),
\]
where $a^{(b)}$ indicates that $a$ is repeated $b$ times. Then,
\[
\begin{array}{l}
 \textbf{x}^T (Tr(G') - Tr(G))\textbf{x }= (n_1  - 1)(n + n_1  - 3)x_1^2  + (n_r  + 1)(n + n_r  - 1)x_r^2 - (n_1  - 1) \\
\;\;\; \;\;\;\;\;\;\;\;\;\;\;\;\; \;\;\;\;\;\;\;\;\;\;\;\;\;\;\;\;\;\;\;\;\;\;\;\times(n + n_1  - 2)x_1^2  - n_r (n + n_r  - 2)x_r^2  - (n + n_1  - 2)x_r^2  \\
\;\;\; \;\;\;\;\;\;\;\;\;\;\;\;\;
\;\;\;\;\;\;\;\;\;\;\;\;\;\;\;\;\;\;\;\;=(2n_r-n_1+1)x_r^2-(n_1-1)x_1^2.\\
 \end{array}
\]
Since $n_1>n_r+1$, then Claim 1 implied that $\textbf{x}^T (Tr(G') -
Tr(G))\textbf{x}<0.$ Similarly, from Claim 1 again, we have
\[
\begin{array}{l}
 \textbf{x}^T (D(G') - D(G))\textbf{x} = 2(n_1  - 1)x_1 x_r  + 4n_r x_r^2  - 4(n_1  - 1)x_1 x_r  - 2n_r x_r^2  \\
  \;\;\; \;\;\;\;\;\;\;\;\;\;\;\;\;
\;\;\;\;\;\;\;\;\;\;\;\;\;\;\;\;\;= 2x_r (n_r x_r  - (n_1  - 1)x_1 ) < 0. \\
 \end{array}
\]
Thus, from the Rayleigh's principle, we have
\[
\begin{array}{l}
  \rho_{D_\alpha}(G')-\rho_{D_\alpha}(G)\le \textbf{x}^T (D_\alpha  (G') - D_\alpha  (G))\textbf{x} \\
  \;\;\; \;\;\;\;\;\;\;\;\;\;\;\;\;
\;\;\;\;\;\;\;\;\;\;\;\;= \alpha x^T (Tr(G') - Tr(G))x + (1 - \alpha )x^T (D(G') - D(G))x \\
 \;\;\; \;\;\;\;\;\;\;\;\;\;\;\;\;
\;\;\;\;\;\;\;\;\;\;\;\; < 0, \\
 \end{array}
\]
contradicting the minimality argument of the generalized distance
spectral radius of $G$. Hence, $G\cong T_{n,r}$.

\paragraph{Claim 2.} {\it
\begin{equation*}
\rho _{D_\alpha  } (T_{n,r})=\frac{{n + 2d + \alpha n - 3 + \sqrt
{(n(1 - \alpha ) - 1)^2  + 4s(1 - \alpha )(d + 1)} }}{2},
\end{equation*}
where $n=rd+s$, $0\leq s<r$.}\\
\\
\emph{Proof of Claim 2.} Let $V_i\;(i=1,2,\ldots,r)$ be the vertex
partition of $T_{n,r}$ and $\textbf{x}$ be the unit Perron vector of
$T_{n,r}$ corresponding to $\rho_{D_\alpha}(T_{n,r})$. From Lemma
2.2, we have $x_r=x_s$ for any two vertices $v_r,v_s$ in every
subset $V_i$. Then the identity (\ref{3}) implies that, for any
$i=1,2,\ldots,r$,
\[
\rho _{D_\alpha  }(T_{n,r}) x_i  = \alpha (n + n_i  - 2)x_i  + (1 -
\alpha )\left( {\sum\limits_{k = 1,k \ne i}^r {n_k x_k }  + 2(n_i  -
1)x_i } \right),
\]
which is equivalent to
\begin{equation}\label{7}
\frac{{n_i }}{{\rho _{D_\alpha  } (T_{n,r}) - \alpha n - n_i  + 2}}
= \frac{{n_i x_i }}{{(1 - \alpha )\sum\nolimits_{k = 1}^r {n_k x_k }
}}.
\end{equation}
Observe that $T_{n,r}$ has $s$ parts of order $d+1$ and $r-s$ parts
of order $d$. Summing the identity (\ref{7}) for all
$i=1,2,\ldots,r$, we see that $\rho _{D_\alpha  }(T_{n,r})$
satisfies the following equation in $\lambda$:
\[
\frac{{s(d + 1)}}{{\lambda  - \alpha n - d + 1}} + \frac{{(r -
s)d}}{{\lambda  - \alpha n - d + 2}} = \frac{1}{{(1 - \alpha )}}.
\]
By solving this equation, we get
\[
\lambda _ \pm   = \frac{{n + 2d + \alpha n - 3 \pm \sqrt {(n(1 -
\alpha ) - 1)^2  + 4s(1 - \alpha )(d + 1)} }}{2}.
\]
On the other hand, from (\ref{7}), we easily see that $\rho
_{D_\alpha  }(T_{n,r})>\alpha n + n_i -2\geq \alpha n +d-1.$ But, $
\lambda _ -   \le \alpha n + d - 1$ as $n>r\geq\frac{1}{1-\alpha}$.
Therefore, $\rho _{D_\alpha  }(T_{n,r})=\lambda_+$.

Above all, the proof of Theorem 2.4 is completed.
\end{proof}

In particular, for $\alpha=0$ in Theorem 2.4, we get a lower bound
on distance spectral radius of connected graphs in terms of
chromatic number.

\paragraph{Corollary 2.5.}(\cite{Liu2010}) {\it Let $G$ be a connected graph of order $n$
with given chromatic number $r$, where $3 \leq r \leq n-1$. Then
\begin{equation}\label{8}
\rho _{D} (G) \ge \rho _{D } (T_{n,r})=\frac{{n + 2d  - 3 + \sqrt
{(n - 1)^2  + 4s(d + 1)} }}{2},
\end{equation}
where $n=rd+s, 0\leq s<r$. Furthermore, the equality in (\ref{8})
holds if and only if $G\cong T_{n,r}$.}\\

Now, for $\alpha=\frac{1}{2}$ in Theorem 2.4, we easily obtain a
lower bound about the distance signless Laplacian spectral radius of
connected graphs in terms of chromatic number.

\paragraph{Corollary 2.6.}(\cite{Li2014}) {\it Let $G$ be a connected graph of order $n$
with given chromatic number $r$, where $3 \leq r \leq n-1$. Then
\begin{equation}\label{9}
\rho _{D^Q} (G) \ge \rho _{D^Q} (T_{n,r})=\frac{{3n + 4d - 6 + \sqrt
{(n - 2)^2  + 8s(d + 1)} }}{2},
\end{equation}
where $n=rd+s, 0\leq s<r$. Furthermore, the equality in (\ref{9})
holds if and only if $G\cong T_{n,r}$.}

\paragraph{Problem.} {\it Characterize the graphs with minimum
generalized distance spectral radius among all connected graphs of
order $n$ with fixed chromatic number $r\geq 3$ for
$1>\alpha>\frac{1}{r}$.}

\section{Minimum generalized distance spectral radius of the $n$-vertex tree and unicyclic graphs}

In this section, we determine the unique graphs with minimum
generalized distance spectral radius among the $n$-vertex tree and
unicyclic graphs, respectively. We first recall some lemmas, which
will be used in the following section.

The Wiener index \cite{Das2010} of a connected graph $G$, denoted by
$W(G)$, is the sum of distances between all unordered pairs of
vertices in $G$, that is, $W(G) = \sum\nolimits_{1 \le i < j \le n}
{d_{v_iv_j}}=\frac{1}{2}\sum\nolimits_{u \in V(G)} {Tr(u)}.$ Denote
the $n$-vertex star by $S_n$. For $n\geq 3$, let $S_n^+$ denote the
$n$-vertex unicyclic graph obtained by adding an edge to $S_n$.

\paragraph{Lemma 3.1.}(\cite{Cui2019}) {\it Let $G$ be a connected graph on $n$ vertices. Then
 $$  \rho_{D_\alpha} (G) \ge \frac{{2W(G)}}{n}$$
with equality holding if and only if $G$ is transmission regular.}

\paragraph{Lemma 3.2.}(\cite{Xing2014}) (i) {\it Let $G$ be a
tree of order $n\geq 4$ different from $S_n$. Then,
\begin{equation*}
W(G) \ge n^2-n-2>W(S_n)=(n-1)^2;
\end{equation*}}
(ii) {\it Let $G$ be a connected unicyclic graph of order $n\geq 6$
different from $S_n^+$. Then,
\begin{equation*}
W(G) \ge n^2-n-4>W(S_n^+)=n^2-2n.
\end{equation*}}

\paragraph{Lemma 3.3.}(\cite{Cui2019}) {\it Let $G$ be a
connected graph of order $n$. If $1\ge \alpha >\beta\ge 0 $, then
\begin{equation}\label{10}
\rho _{D_\alpha  } (G)\ge \rho _{D_\beta  } (G).
\end{equation}
The equality in (\ref{10}) holds if and only if $G$ is transmission
regular.}

\paragraph{Theorem 3.4.} {\it Let $G$ be a
tree of order $n\geq 4$. If $0\leq\alpha<1$, then,
\[
\rho _{D_\alpha  } (G) \ge \rho _{D_\alpha  } (S_n ) = \frac{{\alpha
n + 2(n - 2) + \sqrt {n^2 (2 - \alpha )^2  + 4(n - 1)(2\alpha  - 3)}
}}{2}
\]
with equality holding if and only if $G\cong S_n$.}

\begin{proof}
Suppose that $G$ is a tree of order $n\geq 4$ different from $S_n$.
Now, from the (i) of Lemma 3.2, we get that $W(G) \ge n^2-n-2$. It
follows from Lemma 3.1 that
\begin{equation}\label{11}
\rho _{D_\alpha  } (G) \ge \frac{{2W(G)}}{n} \ge 2n - 2 -
\frac{4}{n} \ge 2n - 3,
\end{equation}
where the last inequality holds as $n\geq 4$. Observe that the
transmission degree sequences of $S_n$ is
$((2n-3)^{(n-1)},(n-1)^{(1)})$. Then, $\rho _{D_1} (S_n)=2n-3$.
Since $S_n$ is non-transmission regular. Then, Lemma 3.3 implies
that $\rho _{D_\alpha} (S_n)$ is strictly increasing function in
$\alpha$, that is, $\rho _{D_1} (S_n)>\rho _{D_\alpha} (S_n)$ for
$0\leq\alpha<1$. Hence, from (\ref{11}), we have $\rho _{D_\alpha  }
(G)>\rho _{D_\alpha} (S_n)$ for $0\leq\alpha<1$.

Next, we shall give the generalized distance spectral radius $\rho
_{D_\alpha  } (S_n )$. Let $\textbf{x}$ be the unit Perron vector of
$S_n$ corresponding to $\rho_{D_\alpha}(S_n)$. By Lemma 2.2, the
equivalent vertices have same entries in $\textbf{x}$. Thus
$\textbf{x}$ can be shaped as
$$\textbf{x}=(x_1,\underbrace {x_2,\ldots ,x_2}_{n-1}).$$
It follows from (\ref{3}) that
\[
\begin{array}{l}
 \rho _{D_\alpha  } (S_n )x_1  = \alpha (n - 1)x_1  + (1 - \alpha )(n - 1)x_2 ,
 \\\\
 \rho _{D_\alpha  } (S_n )x_2  = \alpha (2n - 3)x_2  + (1 - \alpha )(x_1  + 2(n - 2)x_2 ) .\\
 \end{array}
\]
Eliminating $x_1$ and $x_2$ from above equations, we obtain that
$\rho _{D_\alpha  } (S_n )$ is the largest root of the following
equation in $\rho$:
$$\rho ^2  - (\alpha n + 2n - 4)\rho  + (n - 1)(2\alpha n - 2\alpha -
1) = 0.$$ By solving this equation, we get
\[
\rho _{D_\alpha  } (S_n ) = \frac{{\alpha n + 2(n - 2) + \sqrt {(n^2
(2 - \alpha )^2  + 4(n - 1)(2\alpha  - 3)} }}{2},
\]
completing the proof.
\end{proof}

In particular, for $\alpha=0$ in Theorem 3.4, we get a lower bound
on distance spectral radius of the $n$-vertex trees.

\paragraph{Corollary 3.5.}(\cite{Stevanovic2010}) {\it Let $G$ be a
tree of order $n\geq 4$. Then,
\[
\rho _{D} (G) \geq  \rho _{D} (S_n) = {n -2 + \sqrt {(n-2)^2+(n -
1)} }
\]
with equality holding if and only if $G\cong S_n$.}\\

Similarly, for $\alpha=\frac{1}{2}$ in Theorem 3.4, we arrive at:

\paragraph{Corollary 3.6.}(\cite{Xing2014}) {\it Let $G$ be a
tree of order $n\geq 4$ different from $S_n$. Then,
\[
\rho _{D^Q} (G) > \frac{{5 n -8 + \sqrt {9(n-2)^2 +4(n - 1)} }}{2}.
\]}
\paragraph{Theorem 3.7.} {\it Let $G$ be a
connected unicyclic graph of order $n\geq 6$. Then, there exists
$\alpha_0\in (\frac{1}{2},1)$ such that, for any
 $0\leq\alpha<\alpha_0$,
\begin{equation}\label{12}
\rho _{D_\alpha  } (G) \ge \rho_{D_{\alpha}}(S_n^+),
\end{equation}
where $\rho_{D_{\alpha}}(S_n^+)$ is the largest root of the
following equation in $\rho$:
\[
\begin{array}{l}
 \rho ^3  - (3\alpha n + 2n - \alpha  - 7)\rho ^2  + (2\alpha ^2 n^2  + 6\alpha n^2  - \alpha ^2 n - 17\alpha n + 2\alpha  - 7n + 17)\rho  \\
\;\;\;\;- (4\alpha ^2 n^3  - 10\alpha ^2 n^2  - 8\alpha n^2  + 4\alpha ^2 n + 19\alpha n - 7\alpha  + 3n - 5) = 0. \\
 \end{array}
\]
Furthermore, the equality holds if and only if $G\cong S_n^+$.}

\begin{proof} Suppose that $G$ is a connected unicyclic graph of order $n\geq 6$ different from $S_n^+$.
Now, from the (ii) of Lemma 3.2, we get that $W(G) \ge n^2-n-4$. It
follows from Lemma 3.1 that
\begin{equation}\label{13}
\rho _{D_\alpha  } (G) \ge \frac{{2W(G)}}{n} \ge 2n - 2 -
\frac{8}{n}.
\end{equation}

For $n=6,7$, by a simple calculation, we get the Table 1.
\begin{table}[!htb]
\centering
\begin{tabular}
{p{25pt}|p{45pt}|p{60pt}|p{40pt}}
\hline  \quad $n$ &  $\rho_{D_{1/2}}(S_n^+)$ & $2n-2-8/n$ & $\rho_{D_1}(S_n^+)$ \\
\hline  \quad 6 &  8.3574 &8.6667& 9 \\
\hline \quad 7 &  10.4031 & 10.8571& 11\\
\hline
\end{tabular}\\
\centerline{} \centerline{\small Table 1: The value of
$\rho_{D_{1/2}}(S_n^+)$, $\rho_{D_1}(S_n^+)$ and $2n-2-8/n$ for
$n=6,7.$}
\end{table}
Since every eigenvalue of a matrix is a continuous function of its
elements. Then, from Table 1, there exists some real $\alpha_1\in
(\frac{1}{2},1)$ such that
$\rho_{D_{\alpha_1}}(S_n^+)=2n-2-\frac{8}{n}$ for $n=6$. Similarly,
there exists $\alpha_2\in (\frac{1}{2},1)$ such that
$\rho_{D_{\alpha_2}}(S_n^+)=2n-2-\frac{8}{n}$ for $n=7$. Set
$\alpha_0=\min \{\alpha_1,\alpha_2\}$, Lemma 3.3 implies that
$\rho_{D_{\alpha_0}}(S_n^+)\leq 2n-2-\frac{8}{n}$ for $n=6,7.$
Notice that $S_n^+$ is non-transmission regular. Hence, from Lemma
3.3 again, there exists $\alpha_0\in (\frac{1}{2},1)$ such that
$\rho _{D_\alpha  } (G) \ge
\rho_{D_{\alpha_0}}(S_n^+)>\rho_{D_{\alpha}}(S_n^+)$ for
$0\leq\alpha<\alpha_0$.

Now assume that $n\geq 8$. Then, $2n-2-\frac{8}{n}\geq
\rho_{D_1}(S_n^+)=2n-3$. Since $S_n^+$ is non-transmission regular.
Then, from (\ref{13}) and Lemma 3.3, we obtain $\rho _{D_\alpha  }
(G) \ge
\rho_{D_{1}}(S_n^+)>\rho_{D_{\alpha_0}}(S_n^+)>\rho_{D_{\alpha}}(S_n^+)$
for $0\leq\alpha<\alpha_0<1$.

Finally, notice that the transmission degree sequences of $S_n^+$ is
$((2n-3)^{(n-3)},(2n-4)^{(2)},(n-1)^{(1)})$. Then the quotient
matrix of $D_\alpha(S_n^+)$ can be shaped as the form:
\[
M=\left( {\begin{array}{*{20}c}
   {\alpha (n - 1)} & {2(1 - \alpha )} & {(n - 3)(1 - \alpha )}  \\
   {1 - \alpha } & {\alpha (2n - 4) + (1 - \alpha )} & {2(n - 3)(1 - \alpha )}  \\
   {1 - \alpha } & {4(1 - \alpha )} & {\alpha (2n - 3) + 2(n - 4)(1 - \alpha )}  \\
\end{array}} \right).
\]
The spectral radius of a nonnegative square matrix is same as the
spectral radius of a quotient matrix of it corresponding to an
equitable partition (see \cite{Atik2018}), implying the required
results.
\end{proof}

Remark that, for $\alpha=0$ ($\alpha=\frac{1}{2}$, respectively) in
Theorem 3.7, we easily obtain the following Corollary 3.8 (Corollary
3.9, respectively).

\paragraph{Corollary 3.8.}(\cite{Yu2011}) {\it Let $G$ be a
connected unicyclic graph of order $n\geq 6$. Then, $ \rho _{D} (G)
\geq  \rho _{D} (S_n^+) $ with equality holding if and only if
$G\cong S_n^+$.}

\paragraph{Corollary 3.9.}(\cite{Xing2014}) {\it Let $G$ be a
connected unicyclic graph of order $n\geq 6$ different from $S_n^+$.
Then, $ \rho _{D^Q} (G) > \rho _{D^Q} (S_n^+) $.}

\paragraph{Acknowledgements} This work was supported by the
National Natural Science Foundation of China (Nos. 11801521, 11671053).

\end{document}